\UseRawInputEncoding
\documentclass[final,11pt,a4paper,reqno]{amsart}
\usepackage{stmaryrd}
\usepackage[left]{showlabels}
\usepackage[english]{babel}
\usepackage{pst-grad} 
\usepackage{pst-plot} 
\usepackage{pstricks}
\usepackage{amsmath,amssymb,mathrsfs,amsthm, mathtools}
\usepackage{tikz} 
\usepackage{mathcomp,wasysym}  
\usepackage{graphicx,subfigure}  
\usepackage[all]{xy} \xyoption{arc} \xyoption{color}
\usepackage{epsfig}
\usepackage{cite}
\usepackage[a4paper,
left=2.5cm, right=2.5cm,
top=3cm, bottom=3cm]{geometry}
\usepackage[bookmarks]{hyperref}
\hypersetup{
    colorlinks=true,
    linkcolor=cyan,
    filecolor=magenta,      
    urlcolor=cyan,
    citecolor=blue,
    pdftitle={},
    pdfpagemode=FullScreen,
    }
\usepackage{pgfplots}

\newcommand{\norm}[1]{\left\| #1 \right\|}
\newcommand{\av}[1]{\left| #1 \right|}

\newcommand{\RR}{\mathbb R}

\newcommand{\pa}{\partial}

\normalsize
\normalsize
\newtheorem{theorem}{Theorem}

\theoremstyle{definition}

\pgfplotsset{
    every axis/.append style={
            axis x line=middle,
            axis y line=middle,
            xlabel={$x$},
            ylabel={$y$},
            axis line style={->},
        },
    marya/.style={color=black,thick,mark=none},
    soldot/.style={color=black,only marks,mark=*},
    holdot/.style={color=black,fill=white,only marks,mark=*},
    grid style={dotted,gray},
}

\tikzset{>=stealth}

\title[Well-posedness in cold plasmas]{Well-posedness for an hyperbolic-hyperbolic-elliptic system describing cold plasmas}

\author[D. Alonso-Or\'{a}n]{Diego Alonso-Or\'{a}n}
\email{dalonsoo@ull.edu.es}
\address{Departamento de An\'{a}lisis Matem\'{a}tico y Instituto de Matem\'{a}ticas y Aplicaciones (IMAULL), Universidad de La Laguna C/. Astrof\'{i}sico Francisco S\'{a}nchez s/n, 38200 - La Laguna, Spain.}

\author[R. Granero-Belinch\'{o}n]{Rafael Granero-Belinch\'{o}n}
\email{rafael.granero@unican.es}
\address{Departamento  de  Matem\'aticas,  Estad\'istica  y  Computaci\'on,  Universidad  de Cantabria.  Avda.  Los  Castros  s/n,  Santander,  Spain.}

\begin{document}
\begin{abstract}
In this short note, we provide the well-posedness for an hyperbolic-hyperbolic-elliptic system of PDEs describing the motion of collision free-plasma in magnetic fields. The proof combines a pointwise estimate together with a bootstrap type of argument for the elliptic part of the system.
\end{abstract}



\maketitle

\allowdisplaybreaks

\section{Introduction and main result}
The motion of a cold plasma in a magnetic field consisting of singly-charged particles can be described by the following system of PDEs \cite{berezin1964theory, gardner1960similarity}
\begin{subequations}\label{eq:1}
\begin{align}
n_t+(un)_x&=0,\\
u_t+uu_x+\frac{BB_x}{n}&=0,\\
B-n-\left(\frac{B_x}{n}\right)_x&=0,
\end{align}
\end{subequations}
where $n,u$ and $B$ are the ionic density, the ionic velocity and the magnetic field, respectively.  Moreover, it has also been used as a simplified model to describe the motion of collission-free two fluid model where the electron inertial, charge separation and displacement current are neglected and the Poisson equation (\ref{eq:1}c) is initially satisfied, \cite{berezin1964theory, kakutani1968reductive}. In (\ref{eq:1}) the spatial domain $\Omega$ is either $\Omega=\mathbb{R}$ or $\Omega=\mathbb{S}^1$ (\emph{i.e.} $x\in\mathbb{R}$ or $x\in [-\pi,\pi]$ with periodic boundary conditions) and the time variable satisfies $t\in [0,T]$ for certain $0<T\leq\infty$. The corresponding initial-value problem consists of the system \eqref{eq:1} along with initial conditions
\begin{equation}\label{initialdata}
n(x,0)=n_0(x),\;u(x,0)=u_0(x),
\end{equation}
which are assumed to be smooth enough for the purposes of the work.

System \eqref{eq:1} was introduced by Gardner \& Morikawa \cite{gardner1960similarity}. Furthermore, Gardner \& Morikawa formally showed that the solutions of \eqref{eq:1} converge to solutions of the Korteveg-de Vries equation (see also the paper by Su \& Gardner \cite{su1969korteweg}). Berezin \& Karpman extended this formal limit to the case where the wave propagates at angles of certain size with respect to the magnetic field \cite{berezin1964theory}. Later on, Kakutani, Ono, Taniuti \& Wei \cite{kakutani1968reductive} removed the hypothesis on the angle. This formal KdV limit was recently justified by Pu \& Li \cite{pu2019kdv}. Very recently in \cite{AloDurGra}, by means of a multi-scale expansion (cf. \cite{AloranWaves,CGSW18}),  the authors derived three asymptotic models of \eqref{eq:1} and studied several analytical properties of the models: the existence of conserved quantities, the Hamiltonian structure, the well-posedness and the formation of singularities in finite time. More precisely, for the uni-directional model which resembles the well-known Fornberg-Whitham equation (cf. \cite{Fornberg-Whitham-78}), the authors showed that wave-breaking occurs, that is, the formation of an infinite slope in the solution. In \cite{ChenYang}, a new sufficient condition on the initial data is given which exhibits wave breaking extending the previous work \cite{AloDurGra}. To the best of the author's knowledge, although system \eqref{eq:1} has been introduced more than 50 year's ago, the well-posedness of the system has not been studied elsewhere. The result of this work is to give a positive answer to the previous problem and the main theorem reads as follows
\begin{theorem}\label{thm1}
Let $n_0(x)>0$, $n_0(x)-1\in H^2$ and $u_0(x)\in H^3$. Then, there exists a $T>0$ and unique solution of \eqref{eq:1} such that
$$
(n-1,u)\in C([0,T],H^2\times H^3).
$$
\end{theorem}

\subsection*{Notation} For $1\leq p\leq\infty$, let $L^{p}=L^{p}(\RR)$ be the usual normed space of $L^{p}$-functions on $\RR$ with $||\cdot ||_{p}$ as the associated norm. For $s\in\RR$, the inhomogeneous Sobolev space $H^{s}=H^s(\RR)$ is defined as
	\begin{align*}
		H^s(\RR)\triangleq\left\{f\in L^2(\RR):\|f\|_{H^s(\RR)}^2=\int_\RR(1+\xi^2)^s|\widehat{f}(\xi)|^2<+\infty\right\},
	\end{align*}
with norm 
\[ \norm{f}_{H^s}=\norm{f}_{L^2}^2+\norm{f}_{\dot{H}^s}. \]
Moreover,  throughout the paper $C = C(\cdot)$ will denote a positive constant that may depend on fixed parameters  and  $x \lesssim y$ ($x \gtrsim y$) means that $x\le C y$ ($x\ge C y$) holds for some $C$.

\section{Proof of Theorem \ref{thm1}}\label{sec:proof:thm}
The proof follows the classical a priori estimates approach which combines the derivation of useful a priori energy estimates and the use of a suitable approximation procedure via mollifiers (see for instance \cite{AloDurGra}). First, we write system \eqref{eq:1} in the new variables
$n=1+\eta, B=1+b$. Then system \eqref{eq:1} becomes
\begin{subequations}\label{eq:2}
\begin{align}
\eta_t+(u\eta)_x+u_{x}&=0,\\
u_t+uu_x+\frac{(1+b)b_x}{1+\eta}&=0,\\
b-\eta-\left(\frac{b_x}{1+\eta}\right)_x&=0.
\end{align}
\end{subequations}
We are going to find the appropriate energy estimates for the following energy
\begin{equation}\label{energy}
\mathcal{E}(t)=\norm{\eta(t)}^{2}_{H^{2}}+\norm{u(t)}^{2}_{H^{3}}+\displaystyle\max_{x\in \mathbb{R}}\frac{1}{1+\eta(x,t)}.
\end{equation}
In order to estimate the last term in the energy $\mathcal{E}(t)$, we need to derive a pointwise estimate. To that purpose, following \cite{Cordoba-Cordoba-04} and defining
\begin{align*}
m(t)=\displaystyle\min_{x\in\RR}\eta(x,t)=\eta(\underline{x}_{t},t), \mbox{ for } t>0,
\end{align*}
it is easy to check that $m(t)$ is a Lipschitz functions and one has the following bound
$$
\left|m(t)-m(s)\right|\leq \max_{y,z}|\partial_t \eta(y,z)||t-s|.
$$
From Rademacher's theorem it holds that $m(t)$ is differentiable in $t$ almost everywhere and furthermore
\begin{equation}\label{adg19}
m'(t)=\partial_t \eta(\underline{x}_{t},t) \text{ a.e.}
\end{equation}
Then, using (\ref{eq:1}a) and noticing that $n_{x}(\underline{x}_{t},t)=$ we readily see that
\begin{align}\label{expresion}
   m'(t)=-u_{x}(\underline{x}_{t},t)m(t)-u_{x}(\underline{x}_{t},t)=-u_{x}(\underline{x}_{t},t)(1+m(t))
\end{align}
Moreover, since by assumption $m(0)>-1$ we also have that
\begin{equation}\label{pointwise}
m(t)>-1, \quad \mbox{ for } 0<t\ll 1.
\end{equation} 
We remark that this is not a monotonicity statement relying on a sign condition for $u_{x}(\underline{x}_{t},t)$, but just a small in time argument.
Hence, following the argument in \cite{CGO} and using \eqref{expresion} we find that
\begin{align}\label{pointwise2}
  \frac{d}{dt} \left(\displaystyle\max_{x\in \mathbb{R}}\frac{1}{1+\eta(x,t)}\right)=-\frac{\partial_{t}\eta(\underline{x}_{t},t)}{(1+m(t))^2}=\frac{u_{x}(\underline{x}_{t},t)}{1+m(t)}\leq  C (\mathcal{E}(t))^2.
\end{align}
The lower order $L^{2}$ norm of $\eta$ is bounded by
\begin{align}\label{L2:eta}
\frac{1}{2}\frac{d}{dt}\norm{\eta}_{L^{2}}^{2}&\lesssim \norm{\eta}_{L^{2}}^{2}\norm{u_{x}}_{L^\infty}+\norm{\eta}_{L^{2}}\norm{u_{x}}_{L^{2}}
\end{align}
Similarly, we find that
\begin{equation}\label{L2:u}
    \frac{1}{2}\frac{d}{dt}\norm{u}_{L^{2}}^{2}\lesssim \left( 1+\norm{b}_{L^{\infty}}\right)\norm{\frac{b_{x}}{1+\eta}}_{L^{2}}\norm{u}_{L^{2}}
\end{equation}

Testing equation (\ref{eq:2}a) and  (\ref{eq:2}b) with $\pa_{x}^{4}\eta$ and $\pa_{x}^{6}u$ respectively, and integrating by parts we have that
\begin{align}
\frac{1}{2}\frac{d}{dt}\norm{\pa_{x}^{2}\eta}_{L^{2}}^{2}&\lesssim \norm{\eta}_{H^{2}}^{2}\norm{u}_{H^{3}}+\norm{\eta}_{H^{2}}\norm{u}_{H^{3}} \label{H2:eta}, \\
\frac{1}{2}\frac{d}{dt}\norm{\pa_{x}^{3}u}_{L^{2}}^{2}& \lesssim \norm{u}_{H^{3}}^{3}+\left( 1+\norm{b}_{L^{\infty}}\right)\norm{\frac{b_{x}}{1+\eta}}_{H^{3}}\norm{u}_{H^{3}}. \label{H3:u}
\end{align}

Therefore, combining \eqref{L2:eta}-\eqref{H3:u} and using Sobolev embedding and Young's inequality that 
\begin{align}\label{estimate:energy1}
\frac{1}{2}\frac{d}{dt}\left(\norm{\eta}_{H^{2}}^{2}+\norm{u}_{H^{3}}^{2}\right)&\lesssim \norm{\eta}_{H^{2}}^{3}+\norm{u}_{H^{3}}^{3} + \left( 1+\norm{b}_{H^{1}}\right)^{2}\norm{\frac{b_{x}}{1+\eta}}^{2}_{H^{3}}+\norm{u}_{H^{3}}^{2}. 
\end{align}
Moreover, using (\ref{eq:2}c) we find that
\begin{align*}
\norm{\frac{b_{x}}{1+\eta}}^{2}_{\dot{H}^{3}}=\int_{\mathbb{R}} \av{ \left(\frac{b_{x}}{1+\eta}\right)_{xxx}}^{2} dx=\int_{\mathbb{R}} \left(\frac{b_{x}}{1+\eta}\right)_{xxx}\left(b-\eta\right)_{xx} \ dx\leq \norm{\frac{b_{x}}{1+\eta}}_{\dot{H}^{3}}\left( \norm{\eta}_{H^{2}}+\norm{b}_{H^{2}}\right).
\end{align*}
Therefore, we find that
\[ \norm{\frac{b_{x}}{1+\eta}}_{H^{3}} \leq  \norm{\eta}_{H^{2}}+\norm{b}_{H^{2}}. \]
Plugging the previous estimate in \eqref{estimate:energy1} we infer that
\begin{align}\label{est:faltab}
\frac{1}{2}\frac{d}{dt}\left(\norm{\eta}_{H^{2}}^{2}+\norm{u}_{H^{3}}^{2}\right)&\lesssim \norm{\eta}_{H^{2}}^{3}+\norm{u}_{H^{3}}^{3} + \left( 1+\norm{b}_{H^{1}}\right)^{2}(\norm{\eta}_{H^{2}}+\norm{b}_{H^{2}})+\norm{u}_{H^{3}}^{2} \nonumber \\
&\lesssim 1+ \norm{\eta}_{H^{2}}^{3}+\norm{u}_{H^{3}}^{3} +\norm{b}_{H^{2}}^{3}.
\end{align}
To close the energy estimate, we need to compute $\norm{b}_{H^{2}}^{3}$. To that purpose, we first find using the elliptic equation (\ref{eq:2}c) and integrating by parts that
\[\norm{b}_{L^{2}}^{2}=\int_{\mathbb{R}} \eta b \ dx + \int_{\mathbb{R}} \left(\frac{b_{x}}{1+\eta}\right)_{x}b \ dx=\int_{\mathbb{R}} \eta b \ dx - \int_{\mathbb{R}} \frac{b_{x}^{2}}{1+\eta} \ dx  \]
Therefore, using the pointwise estimate  \eqref{pointwise} we find that the last term 
\[ -\int_{\mathbb{R}} \frac{b_{x}^{2}}{1+\eta} \  dx \leq  0, \]
and hence Young's inequality yields
\begin{equation}\label{estimate0}
    \norm{b}^{2}_{L^{2}}\leq \norm{\eta}^{2}_{L^{2}}
\end{equation} 
To compute the higher-order norm, let us first write
\begin{equation}\label{estimate1}
\norm{b_{x}}^{2}_{L^2}=\int_{\mathbb{R}} \frac{1+\eta}{1+\eta}(b_{x})^2 \ dx=-\int_{\mathbb{R}} \frac{b_{x}}{1+\eta} (1+\eta)_{x}b \ dx-\int_{\mathbb{R}} \left(\frac{b_{x}}{1+\eta}\right)_{x} (1+\eta) \ b \ dx=I_1+I_2.
\end{equation} 
Using H\"olders and Young's inequality we readily see that 
\begin{equation}\label{estimate2}
|I_{1}|\leq \norm{b_{x}}_{L^{2}}\norm{\eta_{x}}_{L^{\infty}}\norm{b}_{L^{2}}\norm{\frac{1}{1+\eta}}_{L^{\infty}} \leq  \frac{1}{2}\norm{b_{x}}_{L^{2}}^{2}+ C\norm{\eta}^{2}_{H^{2}}\norm{b}^{2}_{L^{2}}\norm{\frac{1}{1+\eta}}^{2}_{L^{\infty}}.
\end{equation}
On the other hand, using once again the elliptic equation (\ref{eq:2}c) we find that 
\begin{equation}\label{estimate3}
    I_{2}=\int_{\mathbb{R}}(\eta-b)(1+\eta)b \ dx =\int_{\mathbb{R}} \left(\eta b+\eta^2 b-b^{2}(1+\eta) \right) \ dx \leq \norm{b}_{L^{2}}\norm{\eta}_{L^2}+ \norm{b}_{L^{2}}\norm{\eta}_{L^2}\norm{\eta}_{L^{\infty}}.
\end{equation}
Therefore, collecting \eqref{estimate1} -\eqref{estimate3} we infer that
\begin{equation}\label{estimate4}
\norm{b_{x}}^{2}_{L^2}\leq  \frac{1}{2}\norm{b_{x}}_{L^{2}}^{2}+ C\norm{\eta}^{2}_{H^{2}}\norm{b}^{2}_{L^{2}}\norm{\frac{1}{1+\eta}}^{2}_{L^{\infty}}+\norm{b}_{L^{2}}\norm{\eta}_{L^2}+ \norm{b}_{L^{2}}\norm{\eta}_{L^2}\norm{\eta}_{L^{\infty}}
\end{equation} 
and hence using \eqref{estimate0} we conclude that
\begin{align}\label{estimate5}
\norm{b_{x}}^{2}_{L^2}&\lesssim  
\norm{\eta}^{2}_{H^{2}}\norm{\eta}^{2}_{L^{2}}\norm{\frac{1}{1+\eta}}^{2}_{L^{\infty}}+\norm{\eta}_{L^{2}}^{2}+ \norm{\eta}_{L^{2}}^{2}\norm{\eta}_{L^{\infty}}.
\end{align} 
We iterate the previous idea, to provide an estimate $\norm{b_{xx}}_{L^{2}}$. To that purpose, we write
\begin{equation*}\label{estimate6}
\norm{b_{xx}}^{2}_{L^2}=-\int_{\mathbb{R}}\frac{1+\eta}{1+\eta}b_{xxx} b_{x} \ dx=\int_{\mathbb{R}}(1+\eta)b_{xx} \left(\frac{b_{x}}{1+\eta}\right)_{x} \ dx+\int_{\mathbb{R}}(1+\eta)_{x}b_{xx} \frac{b_{x}}{1+\eta} \ dx=J_1+J_2.
\end{equation*} 
Using the elliptic equation (\ref{eq:2}c), we have that
\begin{align}
    J_{1}=\int_{\mathbb{R}}(1+\eta)b_{xx} (b-\eta) \ dx &\leq \norm{b_{xx}}_{L^2} \norm{(1+\eta)(b-\eta)}_{L^{2}} \nonumber  \\
    &\leq \frac{1}{2\epsilon}\norm{b_{xx}}^{2}_{L^{2}} + C_{\epsilon}\left(1+\norm{\eta}_{H^{2}}^{4}+\norm{b}_{L^{2}}^{4} \right)
\end{align}
where in the second inequality we have used the Sobolev embedding and Young's ineqality.
Similarly, 
\begin{align}
    J_{2} &\leq \norm{b_{xx}}_{L^2} \norm{(1+\eta)_{x}\frac{b_{x}}{1+\eta}}_{L^{2}} 
    \leq \frac{1}{2\epsilon}\norm{b_{xx}}^{2}_{L^{2}} + C_{\epsilon}\left(1+\norm{\eta}_{H^{2}}^{8}+\norm{\frac{1}{1+\eta}}_{L^{\infty}}^{8} + \norm{b_{x}}_{L^{2}}^{8} \right).
\end{align}
Therefore taking $\epsilon\ll 1$ (for instance $\epsilon=1/4$),  we find that 
\begin{equation}\label{estimateforbxx}
\frac{1}{2}\norm{b_{xx}}_L^{2}\leq C\left(1+\norm{\eta}_{H^{2}}^{8}+\norm{\frac{1}{1+\eta}}_{L^{\infty}}^{8} + \norm{b_{x}}_{L^{2}}^{8} \right). 
\end{equation}
Hence, estimate \eqref{estimateforbxx} combined with the previous estimates for $\norm{b_{x}}_{L^2}$ given in \eqref{estimate5} 
and $\norm{b}_{L^2}$ in \eqref{estimate0} we conclude that 
\begin{equation}
    \norm{b}^{3}_{H^{2}}\leq C \left(1+\mathcal{E}(t)\right)^{p},
\end{equation}
for some $C>0$ and $p>2$ large enough. The precise power of $p$ can be computed though it is not essential to provide a local-in-time solution. Hence, plugging the previous estimate into \eqref{est:faltab} and taking into account \eqref{pointwise2} we conclude that
\begin{equation}\label{final:est}
\frac{d}{dt}\mathcal{E}(t)\leq C \left(1+\mathcal{E}(t)\right)^{p}
\end{equation}
for some $C>0$ and $p>2$ large enough which ensures a local time of existence $T^{\star}>0$ such that
\[ \mathcal{E}(t)\leq 4 \mathcal{E}(0), \quad \mbox{ for  } 0\leq t \leq T^{\star}.\]
In order to construct the solution, we first define the approximate problems using mollifiers, which reads
\begin{subequations}\label{eq:regularized}
\begin{align}
\eta^{\epsilon}_t+\mathcal{J}_{\epsilon}(\mathcal{J}_{\epsilon}u\mathcal{J}_{\epsilon}\eta)_x+\mathcal{J}_{\epsilon}\mathcal{J}_{\epsilon}u^{\epsilon}_{x}&=0,\\
u^{\epsilon}_t+\mathcal{J}_{\epsilon}\left(\mathcal{J}_{\epsilon}u\mathcal{J}_{\epsilon}u_x\right)+\frac{(1+b)b_x}{1+\eta^{\epsilon}}&=0,\\
b-\eta^{\epsilon}-\left(\frac{b_x}{1+\eta^{\epsilon}}\right)_x&=0.
\end{align}
\end{subequations}
Repeating the previous estimates we find a time of existence $T^{\star}>0$ for the sequence of regularized problems. Using compactness arguments and passing to the limit we conclude the proof of existence. The time continuity for the solution is obtained by classical arguments. On the one hand, the differential equation \eqref{final:est} gives the strong right continuity at $t=0$. Using the change of variables $\hat{t}=-t$, we get the strong left continuity at $t=0$, which combined show the continuity in time of the solution. 
\section*{Acknowledgments}
\noindent  D.A-O is supported by the Spanish MINECO through Juan de la Cierva fellowship FJC2020-046032-I.  R.G-B is supported by the project ``Mathematical Analysis of Fluids
and Applications" Grant PID2019-109348GA-I00 funded by MCIN/AEI/ 10.13039/501100011033 and
acronym ``MAFyA". This publication is part of the project PID2019-109348GA-I00 funded by MCIN/ AEI
/10.13039/501100011033. This publication is also supported by a 2021 Leonardo Grant for Researchers
and Cultural Creators, BBVA Foundation. The BBVA Foundation accepts no responsibility for the opinions, statements, and contents included in the project and/or the results thereof, which are entirely the responsibility of the authors. D.A-O and R. G-B are also supported by the project ``An\'alisis Matem\'atico Aplicado y Ecuaciones Diferenciales" Grant PID2022-141187NB-I00 funded by MCIN/ AEIand acronym ``AMAED".
\begin{footnotesize}

\end{footnotesize}
\vspace{2cm}

\end{document}